\title{Waiter--Client and Client--Waiter colourability games on a $k$--uniform hypergraph and the $k$--SAT game}

\author{
Wei En Tan
\thanks{School of Mathematics, University of Birmingham, Edgbaston, Birmingham B15 2TT, United
Kingdom. Email: WET916@bham.ac.uk.} 
}

\documentclass[12pt]{article}
\usepackage{amsmath,amsthm,amssymb,latexsym,color,epsfig,a4,enumerate}

\newtheorem{theorem}{Theorem} [section]
\newtheorem{lemma}[theorem]{Lemma}

\newtheorem{remark}[theorem]{Remark}

\newtheorem{conjecture}[theorem]{Conjecture}

\begin{document}

\maketitle

\begin{abstract}
Waiter--Client and Client--Waiter games are two--player, perfect information games, with no chance moves, played on a finite set (board) with special subsets known as the winning sets. Each round of the biased $(1:q)$ game begins with Waiter offering $q+1$ previously unclaimed elements of the board to Client, who claims one. The $q$ elements remaining are then claimed by Waiter. If Client fully claims a winning set by the time all board elements have been offered, he wins in the Client--Waiter game and loses in the Waiter--Client game. We give an estimate for the threshold bias of the $(1:q)$ Waiter--Client and Client--Waiter versions of two different games: the non--2--colourability game, played on the complete $k$--uniform hypergraph, and the $k$--SAT game. In particular, we show that the unique value of $q$ at which the winner of the Client--Waiter version of the non--2--colourability game changes is $\frac{1}{n}\binom{n}{k}2^{-k(1+o_k(1))}$ and, for the Waiter--Client version, the corresponding value of $q$ is $\frac{1}{n}\binom{n}{k}2^{\Theta_k(k)}$. Additionally, we show that the threshold bias for the Waiter--Client and Client--Waiter versions of the $k$--SAT game is $\frac{1}{n}\binom{n}{k}$ up to a factor that is exponential and polynomial in $k$ respectively. This shows that these games exhibit the ``probabilistic intuition''.
\end{abstract}

\section{Introduction}\label{intro}
In this paper, we estimate the value of the threshold bias for two biased $(1:q)$ Waiter--Client and Client--Waiter games: the non--2--colourability game, played on the complete $n$--vertex $k$--uniform hypergraph $K_n^{(k)}$, and the $k$--SAT game. These belong to a wider class of games known as ``positional games''.

A positional game is a two--player perfect information game where each player takes turns to claim previously unclaimed (free) elements of a set (board) $X$ until all members of $X$ have been claimed. At this point, the game ends. The winner is determined by the winning criteria of the specific type of positional game in play. Such criteria are defined by a set $\mathcal{F}\subseteq 2^X$ of so--called ``winning sets'' which are known to both players before the game begins. A game with board $X$ and set $\mathcal{F}$ of winning sets is often denoted by the pair $(X,\mathcal{F})$. Popular examples of positional games include Tic--Tac--Toe and Hex \cite{Gale}.

Since the influential papers of Hales and Jewett \cite{HJ}, Lehman \cite{Lehman} and Erd\H{o}s and Selfridge \cite{ES}, positional games have become a widely researched area of combinatorics, and developments in this field have made a significant impact in spheres such as computer science, with regards to the de-randomisation of randomised algorithms (see \emph{e.g.} \cite{ES}). In the literature, the board for these games is most commonly the edge set of a graph or hypergraph. For an extensive survey on positional games, the interested reader may refer to the monographs of Beck \cite{TTT} and Hefetz, Krivelevich, Stojakovi\'{c} and Szab\'{o} \cite{HKSS}.

In biased $(1:q)$ Waiter--Client and Client--Waiter games, where $q$ is a positive integer, the two players, Waiter and Client, play in the following way. At the beginning of each round of the $(1:q)$ Waiter--Client game $(X,\mathcal{F})$, Waiter offers exactly $q+1$ free elements of $X$ to Client. Client claims one of these, and the remaining $q$ elements are then claimed by Waiter. If, in the last round, only $1\leqslant r<q+1$ free elements remain, Waiter claims all of them. Waiter wins the game if he can force Client to fully claim a winning set in $\mathcal{F}$. Otherwise, Client wins. In the $(1:q)$ Client--Waiter game $(X,\mathcal{F})$, each round begins with Waiter offering $1\leqslant t\leqslant q+1$ free elements of $X$ to Client. Client then claims one of these elements, and the remainder of the offering (if any) is claimed by Waiter. In this game, Client wins if he can fully claim a winning set in $\mathcal{F}$, otherwise Waiter wins.

Since these games are finite, perfect information, two--player games with no chance moves and no possibility of a draw, a classical result from Game Theory guarantees a winning strategy (\emph{i.e.} a strategy that, if followed, ensures a win regardless of how the opponent plays). Also, both games are bias monotone in Waiter's bias $q$. This means that, if Client has a winning strategy for a $(1:q)$ Waiter--Client game then Client also has a winning strategy for the same game with bias $(1:q+1)$. The analogous implication is true when Waiter has a strategy to win a $(1:q)$ Client--Waiter game. Thus, for each $(1:q)$ Waiter--Client or Client--Waiter game $(X,\mathcal{F})$, there exists a unique value of $q$ at which the winner of the game changes. This is known as the \emph{threshold bias} of the game.

We give bounds on the threshold bias of two specific types of Waiter--Client and Client--Waiter games in this paper. The first game of interest is the non--2--colourability game $(E(K_n^{(k)}),\mathcal{NC}_2)$ played on the edge set $E(K_n^{(k)})$ of the complete $n$--vertex $k$--uniform hypergraph $K_n^{(k)}$, for some positive integer $k$. In this, the set $\mathcal{NC}_2$ of winning sets is defined to be $$\mathcal{NC}_2=\{F\subseteq E(K_n^{(k)}):\chi(F)>2\},$$
where $\chi$ denotes the weak chromatic number. The second game we consider is the $k$--SAT game $(\mathcal{C}^{(k)}_n,\mathcal{F}_{SAT})$, where $k$ is a positive integer. This is played on the set $\mathcal{C}^{(k)}_n$ of all $2^k\binom{n}{k}$ possible \emph{$k$--clauses}, where each $k$--clause is the disjunction of exactly $k$ non--complementary literals taken from $n$ fixed boolean variables $x_1,\ldots,x_n$. By \emph{literal}, we mean a boolean variable $x_i$ or its negation $\neg x_i$. The set $\mathcal{F}_{SAT}$ of winning sets is defined to be
$$\mathcal{F}_{SAT}=\{\phi:\exists \mathcal{S}\subseteq\mathcal{C}^{(k)}_n\text{ s.t. }\phi=\bigwedge\mathcal{S}\text{ and }\phi\text{ is not satisfiable}\},$$
where $\bigwedge\mathcal{S}$ denotes the conjunction of all $k$--clauses in $\mathcal{S}$. 

These games are interesting from a computational point of view, since the complexity of deciding the 2--colourability or the satisfiability of a random instance is not known. Additionally, since extensive work (see \emph{e.g.} \cite{Coja-Oghlan2012,Ding2015}) has gone into understanding 2--colourable hypergraphs and satisfiable $k$--CNF boolean formulae in the random setting, they are also perfect candidates for developing our understanding of an interesting heuristic known as the \emph{probabilistic intuition}. This suggests that the threshold bias of a game, in which both players play optimally, is roughly the same as it would be, with high probability (\emph{i.e.} with probability tending to 1 as the size of the board tends to infinity), if both players play the game randomly instead. Investigating this notion serves as the main motivation for the research in this paper and we discuss it in greater detail in Section \ref{probint}. The results that follow show that this phenomenon occurs for our chosen games.

\subsection{The Results}\label{results}
\subsubsection{The non--2--colourability game}

We prove that the threshold bias for the $(1:q)$ Waiter--Client and Client--Waiter versions of $(E(K_n^{(k)}),\mathcal{NC}_2)$ is $\frac{1}{n}\binom{n}{k}2^{\Theta_k(k)}$ and $\frac{1}{n}\binom{n}{k}2^{-k(1+o_k(1))}$ respectively. In Section \ref{probint}, we will see that these match the probabilistic intuition.

\begin{theorem}\label{th::WC2col_kgraph}
Let $k$, $q$ and $n$ be positive integers, with $n$ sufficiently large and $k\geqslant 2$ fixed, and consider the $(1:q)$ Waiter--Client non--2--colourability game played on the edge set of the complete $k$--uniform hypergraph $K_n^{(k)}$ on $n$ vertices. If $q\leqslant\binom{\lceil n/2\rceil}{k}\frac{\ln 2}{2((1+\ln 2)n+\ln 2)}$, then Waiter can force Client to build a non--2--colourable hypergraph. Also, if $q\geqslant2^{k/2}e^{k/2+1}k\binom{n}{k}/n$, then Client can keep his hypergraph 2--colourable throughout the game.
\end{theorem}

\begin{theorem}\label{th::CW2col_kgraph}
Let $k$, $q$ and $n$ be positive integers, with $n$ sufficiently large and $k\geqslant 2$ fixed, and consider the $(1:q)$ Client--Waiter non--2--colourability game played on the edge set of the complete $k$--uniform hypergraph $K_n^{(k)}$ on $n$ vertices. If $q\leqslant\binom{\lceil n/2\rceil}{k}\frac{\ln 2}{(1+\ln 2)n}$, then Client can build a non--2--colourable hypergraph. However, when $q\geqslant k^32^{-k+5}\binom{n}{k}/n$, Waiter can ensure that Client has a 2--colourable hypergraph at the end of the game.
\end{theorem}

Thus, for the Waiter--Client version, we have a gap of $(1+o(1))(1+1/\ln 2)k\cdot 2^{3k/2+1}e^{k/2+1}$ between the upper and lower bounds of $q$ and, for the Client--Waiter version, we have a gap of $(1+o(1))(1+1/\ln 2) 2^5k^3$.

\begin{remark}\label{generalr_col}
Our proofs of Theorems \ref{th::WC2col_kgraph} and \ref{th::CW2col_kgraph} generalise easily to the Waiter--Client and Client--Waiter non--$r$--colourability game $(E(K_n^{(k)}),\mathcal{NC}_r)$, for any fixed $r,k\geqslant 2$, where $$\mathcal{NC}_r=\{F\subseteq E(K_n^{(k)}):\chi(F)>r\}.$$ In particular, their generalisation shows that the threshold bias for the $(1:q)$ Waiter--Client and Client--Waiter versions of $(E(K_n^{(k)}),\mathcal{NC}_r)$ is $\frac{1}{n}\binom{n}{k}r^{\Theta_k(k)}$ and $\frac{1}{n}\binom{n}{k}r^{-k(1+o_k(1))}$ respectively. In Section \ref{probint}, we will see that these threshold biases match the probabilistic intuition when $r$ is sufficiently large. For the sake of clarity and simplicity of presentation, and since the generalisation to $r$ colours is straightforward, we only include our proofs for the case $r=2$ in this paper.
\end{remark}

\subsubsection{The $k$--SAT game}

We prove that the threshold bias for the $(1:q)$ Waiter--Client and Client--Waiter versions of $(\mathcal{C}^{(k)}_n,\mathcal{F}_{SAT})$ is $\frac{1}{n}\binom{n}{k}$ up to a factor that is exponential and polynomial in $k$ respectively. This approximately matches the probabilistic intuition, as we shall see in Section \ref{probint}.

\begin{theorem}\label{th::WCkSAT}
Let $k$, $q$ and $n$ be positive integers, with $n$ sufficiently large and $k\geqslant 2$ fixed, and consider the $(1:q)$ Waiter--Client $k$--SAT game played on $\mathcal{C}^{(k)}_n$. When $q\leqslant\binom{n}{k}/(2n)$, Waiter can ensure that the conjunction of all $k$--clauses claimed by Client by the end of the game is not satisfiable. However, when $q\geqslant 2^{3k/2}e^{k/2+1}k\binom{n}{k}/n$, Client can ensure that the conjunction of all $k$--clauses he claims remains satisfiable throughout the game.
\end{theorem}

\begin{theorem}\label{th::CWkSAT}
Let $k$, $q$ and $n$ be positive integers, with $n$ sufficiently large and $k\geqslant 2$ fixed, and consider the $(1:q)$ Client--Waiter $k$--SAT game played on $\mathcal{C}^{(k)}_n$. When $q<\binom{n}{k}/n$, Client can ensure that the conjunction of all $k$--clauses he claims by the end of the game is not satisfiable. However, when $q\geqslant 16k^3\binom{n}{k}/n$, Waiter can ensure that the conjunction of all $k$--clauses claimed by Client is satisfiable throughout the game.
\end{theorem}

Thus, for the Waiter--Client version, we have a gap of $k\cdot 2^{3k/2+1}e^{k/2+1}$ between the upper and lower bounds of $q$ and, for the Client--Waiter version, we have a gap of $16k^3$.

\subsection{The Probabilistic Intuition}\label{probint}

As mentioned previously, the main motivation behind our results is to investigate a heuristic known as the \emph{probabilistic intuition}. This was first employed by Chv\'{a}tal and Erd\H{o}s in \cite{CE} and we illustrate it as follows. Consider a $(1:q)$ Waiter--Client game $(X,\mathcal{F})$ and suppose both players play randomly. Then, at the end of the game, Client has a set of $\lfloor |X|/(q+1)\rfloor$ random elements of the board. If, with high probability (\emph{whp}) this set contains a winning set, then the probabilistic intuition predicts that Waiter has a winning strategy for the $(1:q)$ game. On the other hand, it predicts that Client has a winning strategy for the $(1:q)$ game if \emph{whp} his set of elements at the end of the game does \emph{not} contain a winning set. In particular, suppose that the function $m_{\mathcal{F}}=m_{\mathcal{F}}(|X|)$ is such that, for every $m\gg m_{\mathcal{F}}$ (\emph{i.e.} $\lim_{|X|\rightarrow\infty}m_{\mathcal{F}}/m=0$) any set of $m$ random elements of $X$ contains some $\mathcal{A}\in\mathcal{F}$ \emph{whp} and, for every $m\ll m_{\mathcal{F}}$ any set of $m$ random elements of $X$ does \emph{not} contain any $A\in\mathcal{F}$ \emph{whp}. We refer to $m_{\mathcal{F}}$ as the \emph{threshold} for the property of containing a winning set in $\mathcal{F}$. If the threshold bias $q_{\mathcal{F}}$ of the game $(X,\mathcal{F})$ has the same order of magnitude as $|X|/m_{\mathcal{F}}$, we say that the game exhibits the probabilistic intuition. Surprisingly, research into Waiter--Client and Client--Waiter games played on graphs has found the probabilistic intuition exhibited in many such games whose winning sets are defined by various graph properties (see \emph{e.g.} \cite{HKTpcm, BHKL, BHL, ctcs}). Our results give more examples of Waiter--Client and Client--Waiter games for which this is true.

Random play by both Waiter and Client in the $(1:q)$ non--2--colourability game on $K_n^{(k)}$ leaves Client with a random $n$--vertex $k$--uniform hypergraph at the end of the game. Alon and Spencer \cite{AS} were the first to bound the values $c$ for which the random $n$--vertex $k$--uniform hypergraph with $m=cn$ edges, denoted by $\mathcal{H}_k(n,m)$, is 2--colourable. They showed that the threshold $c_{2,k}$ for $c$ (although only conjectured to exist) satisfies $\tilde{c}\cdot2^k/k^2<c_{2,k}<2^{k-1}\ln2-\ln2/2$ for some small constant $\tilde{c}>0$. The gap of order $k^2$ was reduced by Achlioptas, Kim, Krivelevich and Tetali to order $k$ in \cite{AKKT}. Subsequently, Achlioptas and Moore \cite{AM2002a} improved this further by showing that $c_{2,k}\geqslant2^{k-1}\ln2-\ln2/2-(1+\varepsilon)/2$ for any $\varepsilon>0$ and $k$ sufficiently large. Further still, Coja--Oghlan and Zdeborov\'{a} \cite{Coja-Oghlan2012} improved this lower bound by an additive $(1-\ln 2)/2$. Hence, the edge threshold for the 2--colourability of a random $n$--vertex $k$--uniform hypergraph is $c_{2,k}n=2^{k(1+o_k(1))}n$. Therefore, the probabilistic intuition predicts that the threshold bias for the $(1:q)$ non--2--colourability game $(E(K_n^{(k)}),\mathcal{NC}_2)$ is $\frac{1}{n}\binom{n}{k}2^{-k(1+o_k(1))}$ which matches the threshold bias (up to the error term in the exponent) given by Theorems \ref{th::WC2col_kgraph} and \ref{th::CW2col_kgraph}. Thus, both the $(1:q)$ Waiter--Client and the $(1:q)$ Client--Waiter versions of the non--2--colourability game exhibit the probabilistic intuition.

Such a conclusion also holds in the more general setting, concerning the $r$--colourability of $\mathcal{H}_k(n,cn)$ where $r\geqslant 2$. This is demonstrated by the following bounds on the corresponding threshold $c_{r,k}$ for $c$ described in the literature. By generalising a result of Achlioptas and Naor \cite{AN2005} on $r$--colouring a random graph (2--uniform hypergraph), Dyer, Frieze and Greenhill \cite{DFG2015} proved that $(r-1)^{k-1}\ln(r-1)\leqslant c_{r,k} \leqslant (r^{k-1}-1/2)\ln r.$ The lower bound was subsequently improved by Ayre, Coja--Oghlan and Greenhill \cite{ACG2015} to $(r^{k-1}-1/2)\ln r-\ln 2-1.01\ln r/r$ for sufficiently large $r$. Thus, for such $r$, the edge threshold for the $r$--colourability of a random $n$--vertex $k$--uniform hypergraph is $c_{r,k}n=r^{k(1+o_k(1))}n$. Therefore, the probabilistic intuition predicts that the threshold bias for the $(1:q)$ non--$r$--colourability game $(E(K_n^{(k)}),\mathcal{NC}_r)$ is $\frac{1}{n}\binom{n}{k}r^{-k(1+o_k(1))}$ when $r$ is large, which matches the threshold biases mentioned in Remark \ref{generalr_col} (up to the error term in the exponent).

The $(1:q)$ Waiter--Client and Client--Waiter $k$--SAT games also exhibit the probabilistic intuition. Indeed, Coja--Oghlan and Panagiotou \cite{CojaOghlan2016} found that the threshold for the satisfiability of the conjunction of random $k$--clauses in $\mathcal{C}^{(k)}_n$ is $(2^{k}\ln 2-(1+\ln 2)/2+o_k(1))n$ (see \cite{FP, CR, FS, AM2002b, FW, AM2006, AP} for earlier work). Hence, the probabilistic intuition predicts that the threshold bias for the $(1:q)$ $k$--SAT game $(\mathcal{C}^{(k)}_n,\mathcal{F}_{SAT})$ is $\frac{1}{n}\binom{n}{k}(\ln 2-o_k(1))^{-1}$. This is matched, up to a constant factor, by the lower bound on the threshold bias for the Waiter--Client and Client--Waiter $k$--SAT games given in Theorem \ref{th::WCkSAT} and \ref{th::CWkSAT} respectively. Since the gap between the upper and lower bounds depends only on $k$ (exponentially in the Waiter--Client game and polynomially in the Client--Waiter game), the threshold bias for both versions of the $k$--SAT game has the same order of magnitude as that predicted by the probabilistic intuition.

\subsection{Related Work}

Waiter--Client and Client--Waiter games were first introduced by Beck (see \emph{e.g.} \cite{becksec}) under the names Picker--Chooser and Chooser--Picker. They are interesting for a number of different reasons. Firstly, when Waiter plays randomly in a Waiter--Client (Client--Waiter) game, the game becomes the \emph{avoiding} (\emph{embracing}) Achlioptas process. Secondly, as we have seen, these games often exhibit the probabilistic intuition. Finally, recent findings indicate that these games may be useful in the study of the highly popular Maker--Breaker games (see \emph{e.g.} \cite{becksec, cdm, Bednarska, Knox}).

Since their introduction, much work has gone into finding the threshold bias of Waiter--Client and Client--Waiter games, with winning sets defined by a variety of different graph properties including connectivity, non--planarity and containing a $K_t$--minor (see \emph{e.g.} \cite{HKTpcm, BHKL}). In \cite{HKTpcm}, the $(1:q)$ Waiter--Client and Client--Waiter non--2--colourability games were considered within the more general non--$k$--colourability game, but instead of playing on $E(K_n^{(k)})$, Hefetz, Krivelevich and Tan considered play on the edge set $E(K_n)$ of the complete graph on $n$ vertices. They found that, for both the Waiter--Client and the Client--Waiter versions, the threshold bias has order $\Theta_{n,k}\left(n/(k\log k)\right)$, thereby showing that these games also exhibit the probabilistic intuition. It was by generalising the techniques used in their paper, that we obtained our results stated in Section \ref{results}.

To our knowledge, the $k$--SAT game has not yet been considered in the literature on positional games. However, the Achlioptas process for $k$--SAT has been studied (see \emph{e.g.} \cite{SV, Perkins, Dani2013}).

\section{Preliminaries}

Whenever necessary, we assume that the number $n$ of vertices or boolean variables is sufficiently large. We also omit floor and ceiling signs whenever these are not crucial. 

\subsection{Hypergraph notation}
Let $\mathcal{H}$ be any $k$--uniform hypergraph with vertex set $V(\mathcal{H})$ and edge set $E(\mathcal{H})\subseteq 2^{V(\mathcal{H})}$, where each edge consists of exactly $k$ vertices. For a vertex $v \in V(\mathcal{H})$, let $E_{\mathcal{H}}(v)$ denote the set of edges of $\mathcal{H}$ that contain $v$ and let $d_{\mathcal{H}}(v) = |E_{\mathcal{H}}(v)|$. For a set $S \subseteq V(\mathcal{H})$, let $\mathcal{H}[S]$ denote the hypergraph with vertex set $S$ and edge set $\{e\in E(\mathcal{H}):e\subseteq S\}$. The \emph{maximum degree} of $\mathcal{H}$ is denoted by $\Delta(\mathcal{H}) = \max \{d_{\mathcal{H}}(v) : v \in V(\mathcal{H})\}$.

We say that a set $A \subseteq V(\mathcal{H})$ is \emph{independent} in $\mathcal{H}$ if $\{e\in E(\mathcal{H}):e\subseteq A\} = \emptyset$. The \emph{independence number} of $\mathcal{H}$, denoted by $\alpha(\mathcal{H})$, is the maximum size of an independent set of vertices in $\mathcal{H}$. A subhypergraph $\mathcal{H}'\subseteq\mathcal{H}$ (\emph{i.e.} a hypergraph $\mathcal{H}'$ with $V(\mathcal{H}')\subseteq V(\mathcal{H})$ and $E(\mathcal{H}')\subseteq E(\mathcal{H})$) is a \emph{clique} in $\mathcal{H}$ if every set of $k$ vertices in $V(\mathcal{H}')$ is an edge of $\mathcal{H}'$. We sometimes refer to a clique with $\ell$ vertices as an $\ell$--clique. The \emph{clique number} of $\mathcal{H}$, denoted by $\omega(\mathcal{H})$, is the largest $\ell$ such that $\mathcal{H}$ contains an $\ell$--clique. The \emph{weak chromatic number} of $\mathcal{H}$, denoted by $\chi(\mathcal{H})$, is the smallest integer $k$ for which $V(\mathcal{H})$ can be partitioned into $k$ independent sets. For a set $F\subseteq E(\mathcal{H})$, we abuse notation slightly by using $\chi(F)$ to denote the chromatic number of the hypergraph with vertex set $V(\mathcal{H})$ and edge set $F$.

Let us denote the complete $n$--vertex $k$--uniform hypergraph by $K_n^{(k)}$ (\emph{i.e.} $K_n^{(k)}$ is an $n$--clique). At any given moment in a Waiter--Client or Client--Waiter game, played on $E(K_n^{(k)})$, let $E_W$ and $E_C$ denote the set of edges currently owned by Waiter and Client respectively. We denote the hypergraph with vertex set $V(K_n^{(k)})$ and edge set $E_W$ by $\mathcal{H}_W$, and the hypergraph with vertex set $V(K_n^{(k)})$ and edge set $E_C$ by $\mathcal{H}_C$. Moreover, let $\mathcal{H}_F$ be the hypergraph consisting of all edges of $K_n^{(k)}$ that are free at a given moment.

\subsection{$k$--SAT notation}
Let $V=\{x_1,\ldots,x_n\}$ be a set of $n$ boolean variables. We define $B_i:=\{x_i,\neg x_i\}$ for each $i\in[n]$ and $\mathcal{B}:=\{B_i:i\in [n]\}$. The conjunction of any number of $k$--clauses is called a \emph{$k$--CNF boolean formula}. For a set $\mathcal{A}$ of $k$--clauses, let $S(\mathcal{A})$ denote the set of literals that appear in the $k$--clauses of $\mathcal{A}$. In the case where $\mathcal{A}$ consists of a single $k$--clause $c$, we will abuse notation slightly and write $S(c)$ instead of $S(\{c\})$. For a set $L$ of literals, let $\mathcal{A}[L]$ denote the set of $k$--clauses of $\mathcal{A}$ whose literals lie in $L$. Let $N_{\mathcal{A}}(L)$ denote the set of $k$--clauses in $\mathcal{A}$ that contain at least one literal from $L$ and let $d_{\mathcal{A}}(L)=|N_{\mathcal{A}}(L)|$. Additionally, we use $\bigwedge\mathcal{A}$ to denote the conjunction of all $k$--clauses in $\mathcal{A}$ and, for any set $B\subseteq\cup\mathcal{B}$, we use $\bigvee B$ to denote the disjunction of all literals in $B$.

Let us denote the set of all possible $2^k\binom{n}{k}$ $k$--clauses, each consisting of non--complementary literals taken from $n$ boolean variables $x_1,\ldots,x_n$, by $\mathcal{C}^{(k)}_n$. At any given moment in a Waiter--Client or Client--Waiter game, played on $\mathcal{C}^{(k)}_n$, let $\mathcal{C}_C$ denote the set of all $k$--clauses currently owned by Client and let $\mathcal{C}_F$ denote the set of $k$--clauses that are currently free.

\subsection{Useful Tools}

We will use the following two lemmas which result from a standard application of the Lov\'{a}sz Local Lemma (see \emph{e.g.} Chapter 5 in \cite{ProbMethod}, \cite{Gebauer2009}).

\begin{lemma}\label{degcondition}
Let $\mathcal{H}$ be a $k$--uniform hypergraph with maximum degree $\Delta(\mathcal{H})\leqslant 2^k/(8k)$. Then $\mathcal{H}$ is 2--colourable.
\end{lemma}

\begin{lemma}\label{satcondition}
Let $k\geqslant 2$ be an integer. Any $k$--CNF boolean formula in which no variable appears in more than $2^{k-2}/k$ $k$--clauses is satisfiable.
\end{lemma}


We will also use the following game theoretic tools. The first two apply to the \emph{transversal} game $(X,\mathcal{F}^*)$. For a finite set $X$ and ${\mathcal F}\subseteq 2^X$, the \emph{transversal} family of ${\mathcal F}$ is ${\mathcal F}^* := \{A \subseteq X : A \cap B \neq \emptyset \textrm{ for every } B \in {\mathcal F}\}$. 

\begin{theorem}[\cite{HKT}] \label{th::CwinsCWtrans} 
Let $q$ be a positive integer, let $X$ be a finite set and let $\mathcal{F}$ be a family of subsets of $X$. If
$$
\sum_{A \in \mathcal{F}} \left(\frac{q}{q+1}\right)^{|A|} < 1,
$$
then Client has a winning strategy for the $(1:q)$ Client--Waiter game $(X, \mathcal{F}^*)$.
\end{theorem}

\begin{theorem}[\cite{Bednarska}] \label{th::WaiterBES}
Let $q$ be a positive integer, let $X$ be a finite set and let ${\mathcal F}$ be a family of subsets of $X$. If
$$
\sum_{A \in {\mathcal F}} 2^{-|A|/(2q-1)} < 1/2 \,,
$$
then Waiter has a winning strategy for the $(1 : q)$ Waiter--Client game $(X, {\mathcal F}^*)$.
\end{theorem} 

\begin{theorem}[implicit in~\cite{TTT}] \label{th::ClientBES} 
Let $q$ be a positive integer, let $X$ be a finite set, let ${\mathcal F}$ be a family of (not necessarily distinct) subsets of $X$ and let $\Phi(\mathcal{F}) = \sum_{A \in \mathcal{F}} (q+1)^{-|A|}$. Then, when playing the $(1 : q)$ Waiter--Client game $(X, {\mathcal F})$, Client has a strategy to avoid fully claiming more than $\Phi(\mathcal{F})$ sets in ${\mathcal F}$.  
\end{theorem}

The rest of this paper is organised as follows: in Section \ref{sec::non2col} we prove Theorems \ref{th::WC2col_kgraph} and \ref{th::CW2col_kgraph}. In Section \ref{sec::ksat} we prove Theorems \ref{th::WCkSAT} and \ref{th::CWkSAT}. Finally, in Section \ref{sec::openprob} we present some open problems.

\section{The non--2--colourability game}\label{sec::non2col}

\subsection{The Waiter--Client non--2--colourability game}

\noindent \textbf{Proof of Theorem~\ref{th::WC2col_kgraph}}. Fix $k\geqslant 2$.
\vspace*{-0.1in}
\paragraph*{Waiter's strategy:}
Suppose that $q\leqslant\binom{\lceil n/2\rceil}{k}\frac{\ln 2}{2((1+\ln 2)n+\ln 2)}$. Since $\omega(\mathcal{H}_W)=\alpha(\mathcal{H}_C)$ at the end of the game and, for any hypergraph $\mathcal{H}$ on $n$ vertices, $\chi(\mathcal{H})\alpha(\mathcal{H})\geqslant n$, it suffices to show that Waiter has a strategy to force Client to claim an edge in every $\lceil n/2\rceil$--clique of $K_n^{(k)}$. Let $\mathcal{F}$ denote the set of all $\lceil n/2\rceil$--cliques in $K_n^{(k)}$. Observe that
\begin{align*}
\sum_{A\in\mathcal{F}}2^{-|A|/(2q-1)}\leqslant\binom{n}{\lceil n/2\rceil}2^{-\binom{\lceil n/2\rceil}{k}/(2q)}\leqslant\left(\frac{en}{\lceil n/2\rceil}\right)^{\lceil n/2\rceil}2^{-((1+1/\ln 2)n+1)}<\frac{1}{2}\left(\frac{e}{2^{1/\ln 2}}\right)^n=\frac{1}{2},
\end{align*}
where the second inequality holds by our choice of $q$. Hence, by Theorem \ref{th::WaiterBES}, Waiter can ensure that $\alpha(\mathcal{H}_C)<\lceil n/2\rceil$ to give $\chi(\mathcal{H}_C)>2$ by the end of the game.

\paragraph*{Client's strategy:}
Suppose $q\geqslant 2^{k/2}e^{k/2+1}k\binom{n}{k}/n$. Also, let $\mathcal{F}=\{F:\exists S\subseteq V(K_n^{(k)})\text{ s.t. }S\neq\emptyset, F\subseteq E(K_n^{(k)}[S])\text{ and }|F|=\frac{2}{k}|S|\}$. Observe that,
\begin{align*}
\Phi(\mathcal{F})&=\sum_{A\in\mathcal{F}}(q+1)^{-|A|}\leqslant\sum_{t=k}^n\binom{n}{t}\binom{\binom{t}{k}}{2t/k}q^{-2t/k}\leqslant\sum_{t=k}^n\left[\frac{en}{t}\left(\frac{ek\binom{t}{k}}{2tq}\right)^{2/k}\right]^t\\
&\leqslant\sum_{t=1}^n\left[\frac{en}{t}\left(\frac{et^{k-1}}{2q(k-1)!}\right)^{2/k}\right]^t\leqslant\sum_{t=1}^n\left[\frac{en}{t}\left(\frac{1}{(2e)^{k/2}}\left(\frac{t}{n}\right)^{k-1}\right)^{2/k}\right]^t\\
&=\sum_{t=1}^n\left[\frac{1}{2}\left(\frac{t}{n}\right)^{\frac{2}{k}(k-1)-1}\right]^t<\sum_{t=1}^{\infty}\left[\frac{1}{2}\right]^t=1,
\end{align*}
where the fourth inequality follows from our choice of $q$ and since $n$ is assumed to be sufficiently large. Thus, by Theorem \ref{th::ClientBES}, Client can ensure that, for every $S\subseteq V(\mathcal{H}_C)$, there exists a vertex $v\in S$ with $d_{\mathcal{H}_C[S]}(v)\leqslant 1$ at the end of the game. Thus, in $V(\mathcal{H}_C)$, there exists a vertex $v_1$ contained in at most one edge of $\mathcal{H}_C$. In $V(\mathcal{H}_C)\setminus\{v_1\}$, there exists a vertex $v_2$ contained in at most one edge of $\mathcal{H}_C[V(\mathcal{H}_C)\setminus\{v_1\}]$. Continuing in this way, we obtain an ordering $v_1,\ldots,v_n$ of the vertices of $\mathcal{H}_C$ where, for each $v_i$, there is at most one edge of $\mathcal{H}_C$, consisting of vertices in $\{v_i,\ldots,v_n\}$, that contains $v_i$. Therefore, by colouring the vertices of $\mathcal{H}_C$ greedily from $v_n$ to $v_1$, we may obtain a 2--colouring of $\mathcal{H}_C$.

{\hfill $\Box$ \medskip\\}

\subsection{The Client--Waiter non--2--colourability game}

\noindent \textbf{Proof of Theorem~\ref{th::CW2col_kgraph}}. Fix $k\geqslant 2$.
\vspace*{-0.1in}
\paragraph*{Client's strategy:}
Suppose that $q\leqslant\binom{\lceil n/2\rceil}{k}\frac{\ln 2}{(1+\ln 2)n}$. Since $\omega(\mathcal{H}_W)=\alpha(\mathcal{H}_C)$ at the end of the game and, for any hypergraph $\mathcal{H}$ on $n$ vertices, $\chi(\mathcal{H})\alpha(\mathcal{H})\geqslant n$, it suffices to show that Client has a strategy to claim an edge in every $\lceil n/2\rceil$--clique of $K_n^{(k)}$. Let $\mathcal{F}$ denote the set of all $\lceil n/2\rceil$--cliques in $K_n^{(k)}$. Observe that,
\begin{align*}
\sum_{A \in \mathcal{F}} \left(\frac{q}{q+1}\right)^{|A|}&\leqslant\sum_{A\in\mathcal{F}}2^{-|A|/q}\leqslant\binom{n}{\lceil n/2\rceil}2^{-\binom{\lceil n/2\rceil}{k}/q}\\
&\leqslant\left(\frac{en}{\lceil n/2\rceil}\right)^{\lceil n/2\rceil}2^{-(1+1/\ln 2)n}<\left(\frac{e}{2^{1/\ln 2}}\right)^n=1,
\end{align*}
where our third inequality follows from our choice of $q$. Hence, by Theorem \ref{th::CwinsCWtrans}, Client can ensure that $\omega(\mathcal{H}_W)<\lceil n/2\rceil$, which means that $\chi(\mathcal{H}_C)>2$.

\paragraph*{Waiter's strategy:}
Suppose that $q\geqslant k^32^{-k+5}\binom{n}{k}/n$ and let us first consider when $k\geqslant 7$. By Lemma \ref{degcondition}, it suffices for Waiter to ensure that $\Delta(\mathcal{H}_C)\leqslant2^k/(8k)$ at the end of the game. He achieves this as follows. In the first round, Waiter offers $q+1$ arbitrary free edges. After this, whenever Client claims an edge, say $e$ consisting of vertices $v_1,\ldots,v_k$ ordered arbitrarily, Waiter responds in the following way. Let $F_1\subseteq\{e\in E(\mathcal{H}_F):v_1\in e\}$ with size $|F_1|=\min\{d_{\mathcal{H}_F}(v_1),\lfloor (q+1)/k\rfloor\}$ and, for each $2\leqslant j\leqslant k$, let  $F_j\subseteq\{e\in E(\mathcal{H}_F):v_j\in e\}\setminus\cup\{F_\ell:1\leqslant\ell<j\}$ with size $|F_j|=\min\{|\{e\in E(\mathcal{H}_F):v_j\in e\}\setminus\cup\{F_\ell:1\leqslant\ell<j\}|,\lfloor (q+1)/k\rfloor\}$. Immediately after Client has claimed $e$, Waiter offers all edges in $\cup\{F_i:i\in [k]\}$. (Recall that, in any round of a Client--Waiter game, Waiter may offer less than $q+1$ edges if he desires.) If no free edge touches $e$, Waiter performs his response on an edge that Client claimed earlier on in the game. If no free edges touch any of Client's previously claimed edges, Waiter simply offers $\min\{q+1,|E(\mathcal{H}_F)|\}$ arbitrary free edges. It is clear that, by responding to each of Client's moves in this way, Waiter offers every edge of $K_n^{(k)}$ in the game. We claim that Waiter's strategy ensures $\Delta(\mathcal{H}_C)\leqslant2^k/(8k)$ at the end of the game.

Indeed, let $u\in V(K_n^{(k)})$ be an arbitrary vertex. Every time Client claims an edge containing $u$, Waiter offers at least $\lfloor (q+1)/k\rfloor$ free edges containing $u$, except for perhaps the last time he offers edges at $u$ when there may be less than $\lfloor (q+1)/k\rfloor$ such edges available. Every time Waiter offers edges containing $u$, Client may claim at most one of these. Hence, at the end of the game,
$$d_{\mathcal{H}_C}(u)\leqslant \frac{\binom{n-1}{k-1}}{\lfloor (q+1)/k\rfloor}+1\leqslant\frac{2k\binom{n-1}{k-1}}{q}+1\leqslant\frac{2^k}{8k},$$
where the final inequality follows from our choice of $k$ and $q$.

In the case where $2\leqslant k\leqslant 6$, Waiter still performs the above strategy, ensuring that
$$d_{\mathcal{H}_C}(u)\leqslant \frac{\binom{n-1}{k-1}}{\lfloor (q+1)/k\rfloor}+1<\frac{\binom{n-1}{k-1}k}{q+1-k}+1\leqslant 2,$$
at the end of the game, where the final inequality follows from $k\leqslant 6$, our choice of $q$ and for sufficiently large $n$. Thus, every vertex in Client's hypergraph is contained in at most one edge of $\mathcal{H}_C$. It is clear therefore, that $\mathcal{H}_C$ is 2--colourable.

{\hfill $\Box$ \medskip\\}
\section{The $k$--SAT game}\label{sec::ksat}

\subsection{The Waiter--Client $k$--SAT game}
\noindent \textbf{Proof of Theorem~\ref{th::WCkSAT}}. Fix $k\geqslant 2$.
\vspace*{-0.1in}
\paragraph*{Waiter's strategy:} 
Let $q\leqslant\binom{n}{k}/(2n)$ and let $$\mathcal{F}=\{\mathcal{A}\subseteq \mathcal{C}^{(k)}_n: |S(\mathcal{A})|=n\text{ and }\bigvee B\in\mathcal{A},\forall B\subseteq S(\mathcal{A})\text{ with }|B|=k\}.$$ Note that, since no $k$--clause of $\mathcal{C}^{(k)}_n$ contains a pair of complementary literals, $S(\mathcal{A})$ cannot contain a pair of complementary literals, for each $\mathcal{A}\in\mathcal{F}$. Hence, observe that
\begin{align*}
\sum_{\mathcal{A}\in\mathcal{F}}2^{-|\mathcal{A}|/(2q-1)}<2^{n-\binom{n}{k}/(2q)}=1,
\end{align*}
where the final equality follows from our choice of $q$. So, by Theorem \ref{th::WaiterBES}, Waiter can force Client to claim a $k$--clause in every $\mathcal{A}\in\mathcal{F}$ by the end of the game. Thus, for every partition $(V_1, V_2)$ of the literals in $S(\mathcal{C}^{(k)}_n)$, where no pair of complementary literals lie in the same part, Client owns two $k$--clauses, $c_1$ and $c_2$, satisfying $S(c_i)\subseteq V_i$ for every $i\in[2]$. This means that, regardless of how one assigns the values 0 and 1 to the boolean variables, Client will always have a $k$--clause consisting entirely of 0--valued literals, since every $\{0,1\}$--assignment defines such a partition of $S(\mathcal{C}^{(k)}_n)$. Thus, at the end of the game, $\bigwedge\mathcal{C}_C$ is not satisfiable.

\paragraph*{Client's strategy:} Let $q\geqslant 2^{3k/2}e^{k/2+1}k\binom{n}{k}/n$. Observe that, with
$$\mathcal{F}=\left\{\mathcal{A}:\exists\mathcal{D}\subseteq\mathcal{B}\text{ s.t. }\mathcal{D}\neq\emptyset,\mathcal{A}\subseteq \mathcal{C}^{(k)}_n\left[\bigcup\mathcal{D}\right],\text{ and }|\mathcal{A}|=\frac{2}{k}|\mathcal{D}|\right\},$$
we have
\begin{align*}
\Phi(\mathcal{F})&=\sum_{\mathcal{A}\in\mathcal{F}}(q+1)^{-|\mathcal{A}|}<\sum_{t=k}^n\binom{n}{t}\binom{\binom{t}{k}2^k}{2t/k}q^{-2t/k}\leqslant\sum_{t=k}^n\left[\frac{en}{t}\left(\frac{ek2^k\binom{t}{k}}{2tq}\right)^{2/k}\right]^t\\
&\leqslant\sum_{t=1}^n\left[\frac{en}{t}\left(\frac{et^{k-1}2^{k-1}}{q(k-1)!}\right)^{2/k}\right]^t\leqslant\sum_{t=1}^n\left[\frac{en}{t}\left(\frac{1}{(2e)^{k/2}}\left(\frac{t}{n}\right)^{k-1}\right)^{2/k}\right]^t\\
&=\sum_{t=1}^n\left[\frac{1}{2}\left(\frac{t}{n}\right)^{\frac{2}{k}(k-1)-1}\right]^t<\sum_{t=1}^{\infty}\left[\frac{1}{2}\right]^t=1,
\end{align*}
where the fourth inequality follows from our choice of $q$ and since $n$ is assumed to be sufficiently large. Thus, by Theorem \ref{th::ClientBES}, Client can ensure that, for every $\mathcal{D}\subseteq\mathcal{B}$, there exists some $B\in\mathcal{D}$ such that $|B\cap S(\mathcal{C}_C\left[\bigcup\mathcal{D}\right])|\leqslant 1$ at the end of the game. Consequently, there exists an ordering $B_{i_1},\ldots,B_{i_n}$ of the elements of $\mathcal{B}$ satisfying the following. For every $1\leqslant j< n$, there is a literal $v_{i_j}\in B_{i_j}$ such that every $k$--clause $c\in\mathcal{C}_C\left[\cup\{B_{i_k}:k\geqslant j\}\right]$ satisfies $S(c)\cap B_{i_j}=\{v_{i_j}\}$ or $S(c)\cap B_{i_j}=\emptyset$. Assigning the value, 0 or 1, to the variable $x_{i_j}$ such that $v_{i_j}=1$, for every $j\in [n]$, provides a satisfying truth assignment for $\bigwedge\mathcal{C}_C$.

{\hfill $\Box$ \medskip\\}

\subsection{The Client--Waiter $k$--SAT game}

\noindent \textbf{Proof of Theorem~\ref{th::CWkSAT}}. Fix $k\geqslant 2$.
\vspace*{-0.1in}
\paragraph*{Client's strategy:} Let $q<\binom{n}{k}/n$. With
$$\mathcal{F}=\{\mathcal{A}\subseteq \mathcal{C}^{(k)}_n: |S(\mathcal{A})|=n\text{ and }\bigvee B\in\mathcal{A},\forall B\subseteq S(\mathcal{A})\text{ with }|B|=k\},$$
observe that
\begin{align*}
\sum_{\mathcal{A}\in\mathcal{F}}\left(\frac{q}{q+1}\right)^{|\mathcal{A}|}\leqslant\sum_{\mathcal{A}\in\mathcal{F}}2^{-|\mathcal{A}|/q}\leqslant 2^{n-\binom{n}{k}/q}<1,
\end{align*}
where the final inequality follows from our choice of $q$. Hence, by Theorem \ref{th::CwinsCWtrans}, Client can claim a $k$--clause in every $\mathcal{A}\in\mathcal{F}$ by the end of the game. Thus, for every partition $(V_1, V_2)$ of the literals in $S(\mathcal{C}^{(k)}_n)$, where no pair of complementary literals lie in the same part, Client owns two $k$--clauses, $c_1$ and $c_2$, satisfying $S(c_i)\subseteq V_i$ for every $i\in[2]$. This means that, regardless of how one assigns the values 0 and 1 to the boolean variables, Client will always have a $k$--clause consisting entirely of 0--valued literals, since every $\{0,1\}$--assignment defines such a partition of $S(\mathcal{C}^{(k)}_n)$. Hence, $\bigwedge\mathcal{C}_C$ is not satisfiable at the end of the game.

\paragraph*{Waiter's strategy:} Let $q\geqslant 16k^3\binom{n}{k}/n$ and first consider when $k\geqslant 6$. By Lemma \ref{satcondition}, it suffices to show that Waiter can ensure no variable appears in more than $2^{k-2}/k$ $k$--clauses of $\mathcal{C}_C$ at the end of the game. Equivalently, this means that at most $2^{k-2}/k$ $k$--clauses of $\mathcal{C}_C$ are allowed to contain a literal in $B_i$, for every $i\in[n]$. 

Waiter plays as follows. In the first round, Waiter offers $q+1$ arbitrary free $k$--clauses. After this, Waiter responds to every $k$--clause claimed by Client in the following way. Suppose Client has just claimed the $k$--clause $c$ with $S(c)=\{v_1,\ldots,v_k\}$. Then, for each $i\in[k]$, there exists $j_i\in [n]$ such that $v_i\in B_{j_i}$. Let $\mathcal{D}_1\subseteq N_{\mathcal{C}_F}(B_{j_1})$ with size $|\mathcal{D}_1|=\min\{d_{\mathcal{C}_F}(B_{j_1}),\lfloor(q+1)/k\rfloor\}$, and for each $2\leqslant i\leqslant k$, let $\mathcal{D}_i \subseteq N_{\mathcal{C}_F}(B_{j_i})\setminus\cup\{\mathcal{D}_{\ell}:1\leqslant\ell<i\}$ with size $|\mathcal{D}_i|=\min\{|N_{\mathcal{C}_F}(B_{j_i})\setminus\cup\{\mathcal{D}_{\ell}:1\leqslant\ell<i\}|,\lfloor(q+1)/k\rfloor\}$. Waiter immediately offers all $k$--clauses in $\cup\{\mathcal{D}_i:i\in[k]\}$ (recall that Waiter may offer less than $q+1$ $k$--clauses in a Client--Waiter game). If $\cup\{\mathcal{D}_i:i\in[k]\}=\emptyset$, then Waiter performs the described response on a $k$--clause claimed by Client earlier in the game. If no free $k$--clause contains a literal of any $k$--clause previously claimed by Client, then Waiter simply offers $\min\{q+1,|\mathcal{C}_F|\}$ free $k$--clauses. It is clear that Waiter offers every $k$--clause of $\mathcal{C}^{(k)}_n$ at some point in the game if he plays in this way.

We claim that the described strategy ensures that at most $2^{k-2}/k$ $k$--clauses of $\mathcal{C}_C$ contain a literal in $B_i$, for every $i\in[n]$, at the end of the game. Indeed, fix an arbitrary $i\in [n]$ and consider $B_i$. Each time Client claims a $k$--clause that contains a literal in $B_i$, Waiter offers at least $\lfloor(q+1)/k\rfloor$ free $k$--clauses that also contain a literal in $B_i$, except for perhaps the last time when there may be less than $\lfloor (q+1)/k\rfloor$ such free $k$--clauses. Each time Waiter offers these $k$--clauses, Client claims at most one of them. So at the end of the game, the number of $k$--clauses of $\mathcal{C}_C$ containing a literal of $B_i$ is at most
$$\frac{\binom{n-1}{k-1}2^k}{\lfloor (q+1)/k\rfloor}+1\leqslant \frac{2^{k+1}k\binom{n-1}{k-1}}{q}+1\leqslant \frac{2^{k-2}}{k},$$
where the final inequality follows from our choice of $k$ and $q$.

In the case where $2\leqslant k\leqslant 5$, Waiter still performs the above strategy, ensuring that the number of $k$--clauses of $\mathcal{C}_C$ containing a literal of $B_i$ is at most
$$\frac{\binom{n-1}{k-1}2^k}{\lfloor (q+1)/k\rfloor}+1<\frac{\binom{n-1}{k-1}2^kk}{q+1-k}+1\leqslant 2,$$
at the end of the game, where the final inequality follows from $k\leqslant 5$, our choice of $q$ and for $n$ sufficiently large. Thus, each boolean variable appears in at most one $k$--clause of $\mathcal{C}_C$. It is clear therefore, that $\bigwedge\mathcal{C}_C$ is satisfiable at the end of the game.

{\hfill $\Box$ \medskip\\}

\section{Concluding remarks and open problems}
\label{sec::openprob}
In this paper, we proved that the threshold bias of the $(1:q)$ non--2--colourability game $(E(K_n^{(k)}),\mathcal{NC}_2)$ is $\frac{1}{n}\binom{n}{k}2^{\Theta_k(k)}$ and $\frac{1}{n}\binom{n}{k}2^{-k(1+o_k(1))}$ for the Waiter--Client and Client--Waiter versions respectively. In addition, we showed that the threshold bias for both the $(1:q)$ Waiter--Client and Client--Waiter versions of the $k$--SAT game $(C_n^{(k)},\mathcal{F}_{SAT})$ is $\frac{1}{n}\binom{n}{k}$ up to an exponential and polynomial factor in $k$ respectively. This shows that these games exhibit the probabilistic intuition. However, there is room to improve the bounds on the threshold bias for all four games, especially in the Waiter--Client versions where the gap is exponential in $k$. In particular, we believe that the threshold bias of these games is asymptotically equivalent to that predicted by the probabilistic intuition in the following sense.

\begin{conjecture}
Let the threshold bias for the $(1:q)$ Waiter--Client and Client--Waiter non--2--colourability games $(E(K_n^{(k)}),\mathcal{NC}_2)$ be denoted by $b_{\mathcal{NC}_2}^{WC}$ and $b_{\mathcal{NC}_2}^{CW}$ respectively. Then
$$\lim_{k\rightarrow\infty}\left\{\lim_{n\rightarrow\infty}\frac{1}{n}\binom{n}{k}\frac{(b_{\mathcal{NC}_2}^{WC})^{-1}}{2^{k-1}\ln 2-\ln 2/2}\right\}=\lim_{k\rightarrow\infty}\left\{\lim_{n\rightarrow\infty}\frac{1}{n}\binom{n}{k}\frac{(b_{\mathcal{NC}_2}^{CW})^{-1}}{2^{k-1}\ln 2-\ln 2/2}\right\}=1.$$
\end{conjecture}

\begin{conjecture}
Let the threshold bias for the $(1:q)$ Waiter--Client and Client--Waiter $k$--SAT games $(\mathcal{C}_n^{(k)},\mathcal{F}_{SAT})$ be denoted by $b_{\mathcal{F}_{SAT}}^{WC}$ and $b_{\mathcal{F}_{SAT}}^{CW}$ respectively. Then
$$\lim_{k\rightarrow\infty}\left\{\lim_{n\rightarrow\infty}\frac{1}{n}\binom{n}{k}\frac{(b_{\mathcal{F}_{SAT}}^{WC})^{-1}}{\ln 2}\right\}=\lim_{k\rightarrow\infty}\left\{\lim_{n\rightarrow\infty}\frac{1}{n}\binom{n}{k}\frac{(b_{\mathcal{F}_{SAT}}^{CW})^{-1}}{\ln 2}\right\}=1.$$
\end{conjecture}

Such similarity between the threshold bias and the probabilistic intuition has been observed before in other Waiter--Client and Client--Waiter games. For example, the threshold bias for the Waiter--Client $K_t$--minor game, played on the edge set $E(K_n)$ of the complete graph $K_n$, matches the probabilistic intuition to this degree (see \cite{HKTpcm}). Thus, it would be interesting to see if the same is true for the games studied here.


\bibliographystyle{plain}
\bibliography{WCCW2col_kgraph}

\end{document}